\documentclass{article}

\newtheorem{theorem}{Theorem}
\newtheorem{lemma}{Lemma}
\newtheorem{corollary}{Corollary}

\newcommand{\bt}{\begin{theorem}}
\newcommand{\et}{\end{theorem}}
\newcommand{\bl}{\begin{lemma}}
\newcommand{\el}{\end{lemma}}
\newcommand{\bc}{\begin{corollary}}
\newcommand{\ec}{\end{corollary}}
\newcommand{\pf}{{\bf Proof}.\ }
\newcommand{\bq}{\begin{eqnarray*}}
\newcommand{\eq}{\end{eqnarray*}}
\newcommand{\be}{\begin{eqnarray}}
\newcommand{\ee}{\end{eqnarray}}
\newcommand{\beq}{\begin{equation}}
\newcommand{\eeq}{\end{equation}}
\newcommand{\benum}{\begin{enumerate}}
\newcommand{\eenum}{\end{enumerate}}
\newcommand{\ba}{\begin{array}}
\newcommand{\ea}{\end{array}}

\title{On Erd\H os's elementary method \\
in the asymptotic theory of partitions
}
\author{Melvyn B. Nathanson\thanks{
This work was supported in part by grants from the PSC--CUNY Research
Award Program and the NSA Mathematical Sciences Program.}\\
Department of Mathematics\\
Lehman College (CUNY)\\
Bronx, New York 10468\\
e-mail: nathansn@alpha.lehman.cuny.edu}

\date{}

\begin{document}
\maketitle

\begin{abstract}
Let $m \geq 1$, and let $A$ be the set of all positive integers
that belong to a union of $\ell$ distinct congruence classes modulo $m$.
Let $p_A(n)$ denote the number of partitions of $n$ into parts
belonging to $A$.  It is proved that
\[
\log p_A(n) \sim \pi\sqrt{\frac{2\ell n}{3m}}.
\]
The proof is based on  Erd\H os's elementary method to obtain
the asymptotic formula for the usual partition function $p(n)$.
\end{abstract}

\section{Asymptotic formulas for partition functions}
Let $p(n)$ denote the number of partitions of $n$.
Hardy and Ramanujan~\cite{hard-rama17,hard-rama18}
and, independently, Uspensky~\cite{uspe20}
discovered the asymptotic formula
\beq     \label{logpart:hr}
p(n) \sim \frac{1}{4n\sqrt{3}}\exp{\left(\pi\sqrt{\frac{2n}{3}}\right) }.
\eeq
It follows that
\beq      \label{logpart:hrlog}
\log p(n) \sim \pi\sqrt{\frac{2n}{3}}.
\eeq
The proof of~(\ref{logpart:hr}) uses complex analysis and modular functions.
Erd\H os~\cite{erdo42b} later discovered an elementary proof
of this asymptotic formula;
his argument is complicated, but uses only estimates
for the exponential function and induction from the identity
\[
np(n) = \sum_{ka \leq n} ap(n-ka).
\]
The purpose of this paper is to show that Erd\H os's method,
which is rarely used and almost completely forgotten,
is powerful enough to produce asymptotic estimates
for many partition functions.

Let $A$ be a nonempty set of positive integers, and let $p_A(n)$
the number of partitions of $n$ into parts belonging to the set $A$.
In this paper we consider sets $A$ that are unions of congruence classes.
Let $m \geq 1$, and let $r_1,\ldots, r_{\ell}$ be integers such that
\[
1 \leq r_1 < r_2 < \cdots < r_{\ell} \leq m
\]
and
\beq     \label{logpart:gcd}
(r_1,\ldots,r_{\ell},m) = 1.
\eeq
Let $A$ be the set of all positive integers $a$ such that
$a \equiv r_i\pmod{m}$ for some $i$.
The divisibility condition~(\ref{logpart:gcd}) implies that $p_A(n) \geq 1$
for all sufficiently large integers $n$.
We shall prove that
\[
\log p_A(n) \sim \pi\sqrt{\frac{2\ell n}{3m}}.
\]
This result is not new; it is contained, for example,
in a paper of Meinardus~\cite{mein54a} that is heavily analytic.
We shall prove this result using only Erd\H os's elementary method.

Andrews~\cite[Chapter 6]{andr76} provides
references to asymptotic formulas for various partition functions.
Among the few papers that use Erd\H os's ideas
are Freitag~\cite{frei77}, Grosswald~\cite{gros63},
and Kerawala~\cite{kera69a,kera69b}.
Expositions of Erd\H os's original work
can be found in the books of Grosswald~\cite{gros66},
Hua~\cite{hua82}, and Nathanson~\cite{nath00aa}.

\section{Estimates for sums of exponential functions}
\bl    \label{logpart:lemma1}
If $0 \leq t \leq n$, then
\[
\sqrt{n} - \frac{t}{2\sqrt{n}} - \frac{t^2}{2n^{3/2}}
\leq \sqrt{n-t} \leq \sqrt{n} - \frac{t}{2\sqrt{n}}.
\]
\el

\pf
If $0 \leq x \leq 1$, then
\[
1 - \frac{x}{2} - \frac{x^2}{2} \leq (1-x)^{1/2} \leq 1 - \frac{x}{2}.
\]
The result follows by letting $x = t/n$.

\bl    \label{logpart:lemma2}
If $x > 0$, then
\[
\frac{e^{-x}}{\left( 1 - e^{-x}\right)^2}
< \frac{1}{x^2}.
\]
If $0 < x \leq 1$, then
\[
\frac{e^{-x}}{\left( 1 - e^{-x}\right)^2} > \frac{1}{x^2}-2.
\]
\el

\pf
The power series expansion for $e^{x}$ gives
\bq
e^{x/2}-e^{-x/2}
& = & 2\sum_{k=0}^{\infty} \frac{1}{(2k+1)!}\left(\frac{x}{2}\right)^{2k+1}  \\
& = & x + x^3\sum_{k=1}^{\infty} \frac{x^{2k-2}}{(2k+1)!2^{2k}}.
\eq
If $x > 0$, then
\[
e^{x/2}-e^{-x/2} > x
\]
and so
\[
\frac{e^{-x}}{\left( 1 - e^{-x}\right)^2}
= \frac{1}{\left(e^{x/2}-e^{-x/2}\right)^2}
< \frac{1}{x^2}.
\]
If $0 < x \leq 1$, then
\bq
e^{x/2}-e^{-x/2} 
& < & x + x^3\sum_{k=1}^{\infty} \frac{1}{2^{2k}}  \\
& < & x + x^3  \\
& < & \frac{x}{1-x^2}
\eq
and so
\[
\frac{e^{-x}}{\left( 1 - e^{-x}\right)^2}
= \frac{1}{\left(e^{x/2}-e^{-x/2}\right)^2}
> \left( \frac{1}{x} - x \right)^2
> \frac{1}{x^2}-2.
\]

\bl    \label{logpart:lemma:qseries}
If $0 \leq q < 1$, then
\[
\sum_{v=1}^{\infty} v^3q^v < \frac{6q}{(1-q)^4}.
\]
\el

\pf
Differentiating the power series
\[
\frac{1}{1-q} = \sum_{v=0}^{\infty} q^v,
\]
we obtain
\bq
\frac{1}{(1-q)^2} & = & \sum_{v=0}^{\infty} vq^{v-1},\\
\frac{2}{(1-q)^3} & = & \sum_{v=0}^{\infty} v(v-1)q^{v-2},\\
\frac{6}{(1-q)^4} & = & \sum_{v=0}^{\infty} v(v-1)(v-2)q^{v-3}  \\
& = & \sum_{v=0}^{\infty} (v^3 -3v(v-1) -v)q^{v-3}.
\eq
Therefore,
\bq
\sum_{v=0}^{\infty} v^3 q^{v}
& = & \frac{6q^3}{(1-q)^4} + 3q^2\sum_{v=0}^{\infty} v(v-1)q^{v-2}
      + q\sum_{v=0}^{\infty} vq^{v-1}  \\
& = & \frac{6q^3}{(1-q)^4} + \frac{6q^2}{(1-q)^3} + \frac{q}{(1-q)^2} \\
& = & \frac{q^3+4q^2+q}{(1-q)^4} \\
& < & \frac{6q}{(1-q)^4}.
\eq

\bl   \label{logpart:lemma3}
Let $n$ be a positive integer
and let $c_1$ and $\varepsilon$ be positive real numbers.
Then
\[
\sum_{k=1}^{\infty}
\frac{e^{-\frac{c_1k}{2\sqrt{n}}}}{1-e^{-\frac{c_1k}{2\sqrt{n}}}}
= O\left(  n^{\frac{1}{2}+\varepsilon} \right).
\]
\el

\pf
We apply the Lambert series identity
(Hardy and Wright~\cite[Theorem 310]{hard-wrig79})
\[
\sum_{k=1}^{\infty} \frac{q^k}{1-q^k} =
\sum_{k=1}^{\infty} d(k)q^k,
\]
where $0 < q < 1$ and $d(k)$ is the divisor function.
Let
\[
q = e^{-\frac{c_1k}{2\sqrt{n}}}.
\]
Since
\[
d(k) \ll k^{\varepsilon}
\]
(Hardy and Wright~\cite[Theorem 315]{hard-wrig79}),
and since $e^{-x} \ll x^{-(1+2\varepsilon)}$ for $x \geq c_1/(2\sqrt{n})$,
we have
\[
e^{-\frac{c_1k}{2\sqrt{n}}}
\ll \left( \frac{2\sqrt{n}}{c_1k} \right)^{1+2\varepsilon}
\]
and
\bq
\sum_{k=1}^{\infty}
\frac{e^{-\frac{c_1k}{2\sqrt{n}}}}{1-e^{-\frac{c_1k}{2\sqrt{n}}}}
& = &  \sum_{k=1}^{\infty} d(k) e^{-\frac{c_1k}{2\sqrt{n}}}  \\
& \ll &  \sum_{k=1}^{\infty} k^{\varepsilon}
\left(\frac{2\sqrt{n}}{c_1k}\right)^{1+2\varepsilon}  \\
& = & n^{\frac{1}{2}+\varepsilon} \left(\frac{2}{c_1}\right)^{1+2\varepsilon}
\sum_{k=1}^{\infty}\frac{1}{k^{1+\varepsilon}}  \\
& \ll & n^{\frac{1}{2}+\varepsilon}.
\eq
This completes the proof.

\bl      \label{logpart:bigupper}
Let $m,\ell, r_1,\ldots, r_{\ell}$ be positive integers
such that
\[
1 \leq r_1 < r_2 < \cdots < r_{\ell} \leq m.
\]
Let $A$ be the set of all positive integers $a$ such that
\[
a \equiv r_i\pmod{m}\qquad\mbox{for some $i = 1,\ldots, \ell$.}
\]
Let $\vartheta$ be a real number with
\[
\vartheta > -1.
\]
Let
\[
c_0 = \pi\sqrt{\frac{2\ell}{3m}},
\]
and
\[
c_1 = \left( \sqrt{1+\vartheta} \right)c_0 =
\pi\sqrt{\frac{2(1+\vartheta)\ell}{3m}}.
\]
Then
\[
\sum_{k=1}^{\infty}\sum_{a\in A} ae^{-\frac{c_1ka}{2\sqrt{n}}}
= \frac{n}{1+\vartheta} + O\left( n^{\frac{1}{2}+\varepsilon}  \right)
\]
for every $\varepsilon > 0$.
\el

\pf
Let
\[
q = e^{-\frac{c_1k}{2\sqrt{n}}}.
\]
Then $0 < q < 1$.
For $1 \leq r \leq m$, we have
\bq
\sum_{v=0}^{\infty} (r+mv)e^{-\frac{c_1k(r+mv)}{2\sqrt{n}}}
& = & \sum_{v=0}^{\infty} (r+mv)q^{r+mv}  \\
& = & mq^r\sum_{v=0}^{\infty} vq^{mv} + rq^r\sum_{v=0}^{\infty} q^{mv} \\
& = & \frac{mq^{r+m}}{(1-q^m)^2} + \frac{rq^r}{1-q^m} \\
& = & 
\frac{me^{-\frac{c_1k(r+m)}{2\sqrt{n}}}}
{\left(1-e^{-\frac{c_1km}{2\sqrt{n}}}\right)^2}
+ \frac{re^{-\frac{c_1kr}{2\sqrt{n}}}}{1-e^{-\frac{c_1km}{2\sqrt{n}}}}.
\eq
Therefore,
\bq
\sum_{k=1}^{\infty}\sum_{v=0}^{\infty}
(r+mv)e^{-\frac{c_1k(r+mv)}{2\sqrt{n}}}
& = & 
\sum_{k=1}^{\infty} \frac{me^{-\frac{c_1k(r+m)}{2\sqrt{n}}}}
{\left(1-e^{-\frac{c_1km}{2\sqrt{n}}}\right)^2}
+ \sum_{k=1}^{\infty}
\frac{r e^{-\frac{c_1kr}{2\sqrt{n}}}}{1-e^{-\frac{c_1km}{2\sqrt{n}}}}\\
& = & \sum_{k=1}^{\infty} \frac{me^{-\frac{c_1k(r+m)}{2\sqrt{n}}}}
{\left(1-e^{-\frac{c_1km}{2\sqrt{n}}}\right)^2}
+ O\left(n^{\frac{1}{2}+\varepsilon}  \right),
\eq
since
\[
0 < \sum_{k=1}^{\infty}
\frac{r e^{-\frac{c_1kr}{2\sqrt{n}}}}{1-e^{-\frac{c_1km}{2\sqrt{n}}}}
\leq m\sum_{k=1}^{\infty}
\frac{e^{-\frac{c_1k}{2\sqrt{n}}}}{1-e^{-\frac{c_1k}{2\sqrt{n}}}}
\ll n^{\frac{1}{2}+\varepsilon}
\]
by Lemma~\ref{logpart:lemma3}.

From the definitions of the constants $c_0$ and $c_1$, we obtain
the upper bound
\bq
\sum_{k=1}^{\infty} \frac{me^{-\frac{c_1k(r+m)}{2\sqrt{n}}}}
{\left(1-e^{-\frac{c_1km}{2\sqrt{n}}}\right)^2}
& < & \sum_{k=1}^{\infty} \frac{me^{-\frac{c_1km}{2\sqrt{n}}}}
{\left(1-e^{-\frac{c_1km}{2\sqrt{n}}}\right)^2}   \\
& < & \sum_{k=1}^{\infty} \frac{4n}{c_1^2k^2m}
         \hspace{2cm}\mbox{(by Lemma~\ref{logpart:lemma2})} \\
& = & \frac{4n}{(1+\vartheta)c_0^2m} \sum_{k=1}^{\infty} \frac{1}{k^2}  \\
& = & \frac{4n\pi^2}{6(1+\vartheta)c_0^2m} \\
& = & \frac{n}{(1+\vartheta)\ell}.
\eq
Therefore,
\bq
\sum_{k=1}^{\infty}\sum_{a\in A} ae^{-\frac{c_1ka}{2\sqrt{n}}}
& = & \sum_{i=1}^{\ell}\sum_{k=1}^{\infty}\sum_{v=0}^{\infty}
(r_i+mv)e^{-\frac{c_1k(r_i+mv)}{2\sqrt{n}}}  \\
& < & \sum_{i=1}^{\ell}  \left( \frac{n}{(1+\vartheta)\ell}
+ O\left( n^{\frac{1}{2}+\varepsilon}\right) \right)\\
& = & \frac{n}{1+\vartheta} + O\left( n^{\frac{1}{2}+\varepsilon}\right).
\eq

We compute a lower bound as follows:
\bq
\sum_{k=1}^{\infty} \frac{me^{-\frac{c_1k(r+m)}{2\sqrt{n}}}}
{\left(1-e^{-\frac{c_1km}{2\sqrt{n}}}\right)^2}
& > & m\sum_{k \leq \frac{2\sqrt{n}}{c_1m}} e^{-\frac{c_1kr}{2\sqrt{n}}}
\left( \frac{e^{-\frac{c_1km}{2\sqrt{n}}}}
{\left(1-e^{-\frac{c_1km}{2\sqrt{n}}}\right)^2} \right)  \\
& > & m\sum_{k \leq \frac{2\sqrt{n}}{c_1m}} e^{-\frac{c_1kr}{2\sqrt{n}}}
\left( \frac{4n}{c_1^2k^2m^2} - 2 \right)
                \qquad\mbox{(by Lemma~\ref{logpart:lemma2})} \\
& = & \frac{4n}{c_1^2m}\sum_{k \leq \frac{2\sqrt{n}}{c_1m}}
\frac{e^{-\frac{c_1kr}{2\sqrt{n}}}}{k^2} + O\left(\sqrt{n}\right).
\eq
Since $e^{-x} \geq 1 - x$, we have
\bq
\sum_{k \leq \frac{2\sqrt{n}}{c_1m}} \frac{e^{-\frac{c_1kr}{2\sqrt{n}}}}{k^2}
& \geq & \sum_{k \leq \frac{2\sqrt{n}}{c_1m}} \frac{1}{k^2}
\left( 1 -\frac{c_1kr}{2\sqrt{n}}\right)  \\
& = & \sum_{k \leq \frac{2\sqrt{n}}{c_1m}} \frac{1}{k^2}
- \frac{c_1r}{2\sqrt{n}} \sum_{k \leq \frac{2\sqrt{n}}{c_1m}}\frac{1}{k}  \\
& = & \sum_{k=1}^{\infty} \frac{1}{k^2}
- \sum_{k > \frac{2\sqrt{n}}{c_1m}} \frac{1}{k^2}
- \frac{c_1r}{2\sqrt{n}} \sum_{k \leq \frac{2\sqrt{n}}{c_1m}}\frac{1}{k}  \\
& = & \frac{\pi^2}{6} + O\left(\frac{1}{\sqrt{n}}\right)
+ O\left(\frac{\log n}{\sqrt{n}}\right).
\eq
Therefore,
\bq
\sum_{k=1}^{\infty} \frac{me^{-\frac{c_1k(r+m)}{2\sqrt{n}}}}
{\left(1-e^{-\frac{c_1km}{2\sqrt{n}}}\right)^2}
& > & \frac{4n}{c_1^2m}
\left( \frac{\pi^2}{6} + O\left(\frac{1}{\sqrt{n}}\right)
+ O\left(\frac{\log n}{\sqrt{n}}\right)  \right) + O\left(\sqrt{n}\right) \\
& = & \frac{n}{(1+\vartheta)\ell} + O\left( \sqrt{n}\log n\right),
\eq
and so
\bq
\sum_{k=1}^{\infty}\sum_{a\in A} ae^{-\frac{c_1ka}{2\sqrt{n}}}
& = & \sum_{i=1}^{\ell}\sum_{k=1}^{\infty}\sum_{v=0}^{\infty}
(r_i+mv)e^{-\frac{c_1k(r_i+mv)}{2\sqrt{n}}}  \\
& = & \sum_{i=1}^{\ell}\sum_{k=1}^{\infty}
\frac{me^{-\frac{c_1k(r_i+m)}{2\sqrt{n}}}}
{\left(1-e^{-\frac{c_1km}{2\sqrt{n}}}\right)^2}
+ O\left(n^{\frac{1}{2}+\varepsilon}  \right)  \\                  
& > & \sum_{i=1}^{\ell}  \left( \frac{n}{(1+\vartheta)\ell}
+ O\left( \sqrt{n}\log n\right) \right)
+ O\left(n^{\frac{1}{2}+\varepsilon}  \right)  \\
& = & \frac{n}{1+\vartheta} + O\left(n^{\frac{1}{2}+\varepsilon} \right).
\eq                       
This completes the proof.

The notation $\sum_{ka > n}$ (resp. $\sum_{ka \leq n}$)
means the sum over all positive integers
$k$ and all integers $a\in A$ such that $ka > n$
(resp. $ka \leq n$).

\bl     \label{logpart:S1}
Let $A$ be a set of positive integers, and let $c_1$ and $N_0$
be positive numbers.  For every $n \geq 2N_0$,
\[
0 < \sum_{ka > n-N_0} ae^{-\frac{c_1ka}{2\sqrt{n}}}  \ll \frac{1}{\sqrt{n}}.
\]
\el

\pf
If $n \geq 2N_0$, then $n-N_0 \geq n/2.$
Since
\[
e^{-x} \ll \frac{1}{x^6}  \qquad\mbox{for $x > 0$,}
\]
we have
\bq
\sum_{ka > n-N_0} ae^{-\frac{c_1ka}{2\sqrt{n}}}
& \ll & \sum_{ka > n-N_0} a\left( \frac{2\sqrt{n}}{c_1ka}\right)^{6} \\
& \ll & n^3\sum_{ka > n-N_0} \frac{1}{k^6a^5} \\
& \ll & n^3\sum_{ka > n-N_0} \frac{1}{(ka)^{7/2}k^{5/2}a^{3/2}} \\
& \ll & n^3\sum_{ka > n-N_0} \frac{1}{(n/2)^{7/2}k^{5/2}a^{3/2}} \\
& \ll & \frac{1}{\sqrt{n}}\sum_{k=1}^{\infty} \frac{1}{k^{5/2}}
\sum_{a\in A} \frac{1}{a^{3/2}}\\
& \ll & \frac{1}{\sqrt{n}}.
\eq

\bl     \label{logpart:S2}
Let $A$ be a set of positive integers, and let $c_1$ be a positive number.
Then
\[
0 < \sum_{ka \leq n}k^2a^3 e^{-\frac{c_1ka}{2\sqrt{n}}}  \ll n^2.
\]
\el

\pf
This is a straightforward computation.  We have
\bq
\sum_{ka \leq n}k^2a^3 e^{-\frac{c_1ka}{2\sqrt{n}}}
& \leq & \sum_{k=1}^{n} k^2
\sum_{a \in A} a^3 e^{-\frac{c_1ka}{2\sqrt{n}}} \\
& \leq & \sum_{k=1}^{n} k^2
\sum_{v=1}^{\infty}v^3 e^{-\frac{c_1kv}{2\sqrt{n}}} \\
& \leq & 6 \sum_{k=1}^{n}
\frac{k^2 e^{-\frac{c_1k}{2\sqrt{n}}}}
{\left( 1 - e^{-\frac{c_1k}{2\sqrt{n}}} \right)^4 }
\hspace{2cm}\mbox{(by Lemma~\ref{logpart:lemma:qseries})}  \\
& = & 6 \sum_{k=1}^{n}
\frac{e^{-\frac{c_1k}{2\sqrt{n}}}}
{\left( 1 - e^{-\frac{c_1k}{2\sqrt{n}}} \right)^2 }
\frac{k^2}{\left( 1 - e^{-\frac{c_1k}{2\sqrt{n}}} \right)^2} \\
& < & 6 \sum_{k=1}^{n}
\left( \frac{4n}{c_1^2k^2} \right)
\frac{k^2}{\left( 1 - e^{-\frac{c_1k}{2\sqrt{n}}} \right)^2}
\qquad\mbox{(by Lemma~\ref{logpart:lemma2})}  \\
& \ll & n\sum_{k=1}^{n}
\frac{1}{\left( 1 - e^{-\frac{c_1k}{2\sqrt{n}}} \right)^2}.
\eq
Let
\[
x = \frac{c_1k}{2\sqrt{n}}.
\]
If $1 \leq k \leq \sqrt{n}$, then $0 < x \leq c_1/2$ and
\[
1 - e^{-x} = \int_0^x e^{-t}dt \geq xe^{-x} \geq xe^{-c_1/2}.
\]
It follows that
\[
\left(1-e^{-\frac{c_1k}{2\sqrt{n}}}\right)^2
= \left(1-e^{-x}\right)^2 \geq e^{-c_1}x^2
= \frac{e^{-c_1}c_1^2k^2}{4n},
\]
and so
\[
\sum_{1 \leq k \leq \sqrt{n}}
\frac{1}{\left( 1 - e^{-\frac{c_1k}{2\sqrt{n}}} \right)^2}
\leq \frac{4e^{c_1}n}{c_1^2}
\sum_{1 \leq k \leq \sqrt{n}}\frac{1}{k^2}  \ll n.
\]
If $k > \sqrt{n}$, then
\[
\sum_{\sqrt{n} < k \leq n}
\frac{1}{\left( 1 - e^{-\frac{c_1k}{2\sqrt{n}}} \right)^2}
\leq \sum_{\sqrt{n} < k \leq n}
\frac{1}{\left( 1 - e^{-\frac{c_1}{2}} \right)^2} \ll  n.
\]
Therefore,
\[
\sum_{k=1}^{n} \frac{1}{\left( 1 - e^{-\frac{c_1k}{2\sqrt{n}}} \right)^2} \ll n
\]
and
\[
\sum_{ka \leq n}k^2a^3 e^{-\frac{c_1ka}{2\sqrt{n}}}
\ll n\sum_{k=1}^{n}
\frac{1}{\left( 1 - e^{-\frac{c_1k}{2\sqrt{n}}} \right)^2}
\ll n^2.
\]
This completes the proof.

\section{Upper and lower bounds for $\log p_A(n)$}

We define $p_A(0) = 1$ and $p_A(-n) = 0$ for all $n \geq 1$.
We use $k$ to denote a positive integer, $v$ a nonnegative integer,
and $a$ an element of the set $A$ of congruence classes modulo $m$.
The asymptotic formula for $\log p_A(n)$ will be proved by induction
from the following classical recursion formula.

\bl      \label{logpart:recursionformula}
Let $A$ be a nonempty set of positive integers, and let $p_A(n)$
be the number of partitions of $n$ into parts belonging to $A$.
Then
\[
np_A(n) = \sum_{ka \leq n} ap_A(n-ka).
\]
\el

\pf
We enumerate the partitions of $n$ into parts belonging to $A$ as follows:
\[
n = a_{i,1} + a_{i,2} + \cdots + a_{i,s_i}
\qquad\mbox{for $i = 1,\ldots, p_A(n)$}.
\]
Then
\[
np_A(n) = \sum_{i=1}^{p_A(n)} \sum_{j=1}^{s_i} a_{i,j}
= \sum_{a\in A} aN(a,n),
\]
where $N(a,n)$ is the total number of times that the integer $a$ occurs
in the $p_A(n)$ partitions of $n$.
The number of partitions in which the integer $a$
occurs at least $k$ times is $p_A(n-ka)$, and so the number
of partitions in which the integer $a$ occurs {\em exactly} $k$ times is
\[
p_A(n-ka) - p_A(n-(k+1)a).
\]
Therefore,
\[
N(n,a)
= \sum_{k=1}^{\infty} k \left( p_A(n-ka) - p_A(n-(k+1)a) \right)
= \sum_{k=1}^{\infty} p_A(n-ka),
\]
and so
\[
np_A(n) = \sum_{a\in A} aN(a,n)
= \sum_{a\in A} \sum_{k=1}^{\infty} ap_A(n-ka)
= \sum_{ka \leq n} ap_A(n-ka),
\]                           
since $p_A(n-ka) = 0$ if $ka > n$.
This completes the proof.

\begin{theorem}                  \label{logpart:theorem:upper}
Let $m,\ell, r_1,\ldots, r_{\ell}$ be positive integers
such that
\[
1 \leq r_1 < r_2 < \cdots < r_{\ell} \leq m.
\]
Let $A$ be the set of all positive integers $a$ such that
\[
a \equiv r_i\pmod{m}\qquad\mbox{for some $i = 1,\ldots, \ell$.}
\]
Then
\[
\limsup_{n\rightarrow\infty}
\frac{\log p_A(n)}{\pi\sqrt{\frac{2\ell n}{3m}}} \leq 1.
\]
\end{theorem}

\pf            
Let $0 < \varepsilon < 1/2$,
\[
c_0 = \pi\sqrt{\frac{2\ell}{3m}},
\]
and
\[
c_1 = \left(\sqrt{1+\varepsilon}\right)c_0 =
\pi\sqrt{\frac{2(1+\varepsilon)\ell}{3m}}.
\]

We shall prove that there exists a constant $K = K(\varepsilon)$ such that
\beq      \label{logpart:upperlogineq}
p_A(n) \leq Ke^{c_1\sqrt{n}}
\eeq
for all nonnegative integers $n$.
This implies that
\[
\log p_A(n) \leq \log K + \left(\sqrt{1+\varepsilon}\right)c_0\sqrt{n},
\]
and so
\[
\frac{\log p_A(n)}{c_0\sqrt{n}}
\leq \sqrt{1 + \varepsilon} + \frac{\log K}{c_0\sqrt{n}}
\]
and
\[
\limsup_{n\rightarrow\infty}   \frac{\log p_A(n)}{c_0\sqrt{n}} \leq 1.
\]
Therefore, it suffices to prove~(\ref{logpart:upperlogineq}).

Applying Lemma~\ref{logpart:bigupper} with $\vartheta = \varepsilon$,
we have
\[
\sum_{k=1}^{\infty}\sum_{a \in A} ae^{-\frac{c_1ka}{2\sqrt{n}}}
= \frac{n}{1+\varepsilon} + O\left( n^{\frac{1}{2}+\varepsilon}\right).
\]
There exists a positive integer $N = N(\varepsilon)$ such that
\beq       \label{logpart:uppern}
\sum_{k=1}^{\infty}\sum_{a\in A}ae^{-\frac{c_1ka}{2\sqrt{n}}} < n
\eeq
for all $n \geq N$.
We can choose a number $K = K(\varepsilon)$ so that
the upper bound~(\ref{logpart:upperlogineq}) holds
for all positive integers $n \leq N$.
Let $n > N$ and assume that inequality~(\ref{logpart:upperlogineq})
holds for all integers less than $n$.
Then
\bq
np_A(n)
& = & \sum_{k=1}^{\infty}\sum_{a \in A} ap_A(n-ka)
      \hspace{3cm}\mbox{(by Lemma~\ref{logpart:recursionformula})}\\
& \leq & \sum_{k=1}^{\infty}\sum_{a \in A} aKe^{c_1\sqrt{n-ka}}
          \hspace{3cm}\mbox{(by inequality~(\ref{logpart:upperlogineq}))}\\
& \leq & K\sum_{k=1}^{\infty}\sum_{a \in A}
         ae^{c_1\sqrt{n}-\frac{c_1ka}{2\sqrt{n}}}
         \hspace{3cm}\mbox{(by Lemma~\ref{logpart:lemma1})}\\
& = & Ke^{c_1\sqrt{n}}\sum_{k=1}^{\infty}\sum_{a\in A}
         ae^{-\frac{c_1ka}{2\sqrt{n}}} \\
& < & nKe^{c_1\sqrt{n}}
\hspace{4cm}\mbox{(by inequality~(\ref{logpart:uppern}))}.
\eq
Dividing by $n$, we obtain~(\ref{logpart:upperlogineq}).
This completes the proof.

Note that in Theorem~\ref{logpart:theorem:upper}
we do not assume that $(r_1,\ldots,r_{\ell},m) = 1$.

\begin{theorem}         \label{logpart:theorem:lower}
Let $m, \ell, r_1,\ldots, r_{\ell}$ be positive integers
such that
\[
1 \leq r_1 < r_2 < \cdots < r_{\ell} \leq m,
\]
and
\[
(r_1,\ldots,r_{\ell},m) = 1.
\]
Let $A$ be the set of all positive integers $a$ such that
\[
a \equiv r_i\pmod{m}\qquad\mbox{for some $i = 1,\ldots, \ell$.}
\]
Then
\[
\liminf_{n\rightarrow\infty}
\frac{\log p_A(n)}{\pi\sqrt{\frac{2\ell n}{3m}}} \geq 1.
\]
\end{theorem}

\pf            
Let $0 < \varepsilon < 1/2$,
\[
c_0 = \pi\sqrt{\frac{2\ell}{3m}},
\]
and
\[
c_1 = \left(\sqrt{1-\varepsilon}\right)c_0 =
\pi\sqrt{\frac{2(1-\varepsilon)\ell}{3m}}.
\]
The divisibility condition $(r_1,\ldots,r_{\ell},m) = 1$ implies
that there exists a number $N_0$ such that
$p_A(n) \geq 1$ for all integers $n \geq N_0$.
We shall prove that there exists a positive number $K$ such that
\beq      \label{logpart:lowerlogineq}
p_A(n) \geq Ke^{c_1\sqrt{n}}
\eeq
for all integers $n \geq N_0$.
This implies that
\[
\liminf_{n\rightarrow\infty}
\frac{\log p_A(n)}{\pi\sqrt{\frac{2\ell n}{3m}}} \geq 1.
\]

By Lemma~\ref{logpart:bigupper} (with $\vartheta = -\varepsilon$),
Lemma~\ref{logpart:S1}, and Lemma~\ref{logpart:S2},
there exists a number $N_1 = N_1(\varepsilon) \geq 2N_0$
such that, for all integers $n \geq N_1,$
\bq
\lefteqn{
\sum_{ka \leq n - N_0} ae^{-\frac{c_1ka}{2\sqrt{n}}}
+ \frac{c_1}{2n^{3/2}}\sum_{ka \leq n - N_0} k^2a^3e^{-\frac{c_1ka}{2\sqrt{n}}} }\\
& = & \sum_{k=1}^{\infty}\sum_{a\in A} ae^{-\frac{c_1ka}{2\sqrt{n}}}
- \sum_{ka > n - N_0} ae^{-\frac{c_1ka}{2\sqrt{n}}}
+ \frac{c_1}{2n^{3/2}}\sum_{ka\leq n - N_0} k^2a^3e^{-\frac{c_1ka}{2\sqrt{n}}}\\
& = &  \frac{n}{1-\varepsilon}
+ O\left(n^{\frac{1}{2}+\varepsilon}\right)
+ O\left(n^{-\frac{1}{2}}\right)
+ O\left(n^{\frac{1}{2}}\right)  \\
& > & n.
\eq
We can choose a positive number $K = K(\varepsilon)$
such that $p_A(n)$ satisfies
inequality~(\ref{logpart:lowerlogineq}) for $N_0 \leq n \leq N_1$.

Let $n > N_1$, and suppose that inequality~(\ref{logpart:lowerlogineq})
holds for all integers in the interval $[N_0,n-1]$.
We shall prove by induction that this inequality also holds for $n.$
Note that $n-ka \geq N_0$ if $ka \leq n-N_0$.
By Lemma~\ref{logpart:recursionformula}, we have
\bq
np_A(n)
& = & \sum_{ka \leq n} ap_A(n-ka)  \\
& \geq & \sum_{ka \leq n - N_0} ap_A(n-ka)  \\
& \geq & K\sum_{ka \leq n - N_0} ae^{c_1\sqrt{n-ka}}  \\
& \geq & K\sum_{ka \leq n - N_0} ae^{c_1\left(\sqrt{n}-\frac{ka}{2\sqrt{n}}
-\frac{k^2a^2}{2n^{3/2}}  \right)}
\hspace{3cm}\mbox{(by Lemma~\ref{logpart:lemma1})}\\
& = & Ke^{c_1\sqrt{n}}\sum_{ka \leq n - N_0} ae^{-\frac{c_1ka}{2\sqrt{n}}}
e^{-\frac{c_1k^2a^2}{2n^{3/2}}}   \\
& \geq & Ke^{c_1\sqrt{n}}\sum_{ka \leq n - N_0} ae^{-\frac{c_1ka}{2\sqrt{n}}}
\left(1 - \frac{c_1k^2a^2}{2n^{3/2}}\right)
\hspace{0.7cm}\mbox{(since $e^{-x}\geq 1-x$)}\\
& = & Ke^{c_1\sqrt{n}}
\left(\sum_{ka \leq n - N_0} ae^{-\frac{c_1ka}{2\sqrt{n}}}
- \frac{c_1}{2n^{3/2}}\sum_{ka \leq n - N_0} k^2a^3e^{-\frac{c_1ka}{2\sqrt{n}}}
\right) \\
& > & nKe^{c_1\sqrt{n}}.
\eq
Dividing by $n$, we obtain the lower bound~(\ref{logpart:lowerlogineq}).
This completes the induction.

\begin{theorem}            \label{logpart:theorem:main}
Let $m, \ell, r_1,\ldots, r_{\ell}$ be positive integers
such that
\[
1 \leq r_1 < r_2 < \cdots < r_{\ell} \leq m,
\]
and
\[
(r_1,\ldots,r_{\ell},m) = 1.
\]
Let $A$ be the set of all positive integers $a$ such that
\[
a \equiv r_i\pmod{m}\qquad\mbox{for some $i = 1,\ldots, \ell$.}
\]
Then
\[
\log p_A(n) \sim \pi\sqrt{\frac{2\ell n}{3m}}.
\]
\end{theorem}

\pf
This follows immediately from Theorem~\ref{logpart:theorem:upper}
and Theorem~\ref{logpart:theorem:lower}.

The following well--known result is an immediate consequence
of Theorem~\ref{logpart:theorem:main}.

\begin{theorem}
Let $q(n)$ denote the number of partitions of $n$ into distinct parts.
Then
\[
\log q(n) \sim \pi\sqrt{\frac{n}{3}}.
\]
\end{theorem}

\pf
Let $A$ be the set of positive odd numbers, that is, the set of positive
integers congruent to 1 modulo $2$.  Euler proved that $q(n) = p_A(n)$,
so the result follows immediately from
Theorem~\ref{logpart:theorem:main} with $m = 2$ and $\ell = 1$.

\end{document}